\documentclass[11pt]{llncs}
\usepackage{epsfig}
\usepackage{graphicx}
\usepackage{color}
\usepackage{amsmath,amssymb}    
\usepackage{latexsym}
\usepackage{hhline}
\spnewtheorem{observation}{Observation}{\bfseries}{\itshape}
\spnewtheorem{fact}{Fact}{\bfseries}{\itshape}
  \baselineskip 8 mm

\newcommand{\comp}[1]{\overline{#1}}  

\newcommand{\cub}{\mbox{\textnormal{cub}}}
\newcommand{\boxi}{\mbox{\textnormal{box}}}
\newcommand{\cd}{\mbox{\textnormal{chord}}}
\newcommand{\ceil}[1]{\left\lceil #1 \right\rceil}

\newcommand{\claw}{\psi(G)}

\begin{document}

\title{Boxicity and Cubicity of Asteroidal Triple free graphs}
\author {Diptendu Bhowmick
\thanks {Computer Science and Automation Department, Indian
Institute of Science, Bangalore- 560012 Email:
diptendubhowmick@gmail.com},
L. Sunil Chandran
\thanks{Computer Science and Automation Department, Indian Institute
of Science, Bangalore- 560012 Email: sunil@csa.iisc.ernet.in} }
\date{}
\institute{}
\maketitle

\begin{abstract}
An axis parallel $d$-dimensional box is the Cartesian product $R_1 \times R_2 \times \cdots \times R_d$ where each $R_i$ is a closed interval on the real line. The {\it boxicity} of a graph $G$, denoted as $\boxi(G)$, is the minimum integer $d$ such that $G$ can be represented as the intersection graph of a collection of $d$-dimensional boxes. An axis parallel unit cube in $d$-dimensional space or a $d$-cube is defined as the Cartesian product $R_1 \times R_2 \times \cdots \times R_d$ where each $R_i$ is a closed interval on the real line of the form $[a_i,a_i + 1]$. The {\it cubicity} of $G$, denoted as $\cub(G)$, is the minimum integer $d$ such that $G$ can be represented as the intersection graph of a collection of $d$-cubes.

 An independent set of three vertices is called an asteroidal triple if between each pair in the triple there exists a path which avoids the neighbourhood of the third. A graph is said to be Asteroidal Triple free (AT-free for short) if it does not contain an asteroidal triple. The class of AT-free graphs is a reasonably large one, which properly contains the class of interval graphs, trapezoid graphs, permutation graphs, cocomparability graphs etc. Let $S(m)$ denote a star graph on $m+1$ nodes. We define {\it claw number} of a graph $G$ as the largest positive integer $k$ such that $S(k)$ is an induced subgraph of $G$ and denote it as $\claw$.

Let $G$ be an AT-free graph with chromatic number $\chi(G)$ and claw number $\claw$. In this paper we will show that $\boxi(G) \leq \chi(G)$ and this bound is tight. We also show that $\cub(G) \leq \boxi(G)(\ceil{\log_2 \claw} +2)$ $\leq$ $\chi(G)(\ceil{\log_2 \claw} +2)$. If $G$ is an AT-free graph having girth at least $5$ then $\boxi(G) \leq 2$ and therefore $\cub(G) \leq 2\ceil{\log_2 \claw} +4$.\\

\end{abstract}

\section{Introduction}
Let $\mathcal{F}$ be a family of non-empty sets. An undirected graph $G$ is the intersection graph of $\mathcal{F}$ if there exists a one-one correspondence between the vertices of $G$ and the sets in $\mathcal{F}$ such that two vertices in $G$ are adjacent if and only if the corresponding sets have non-empty intersection. If $\mathcal{F}$ is a family of intervals on the real line, then $G$ is called an {\it interval graph}. An interval graph $G$ is said to be a \emph{unit interval graph} if and only if there is some interval representation of $G$ in which all the intervals are of the same length.\\

\noindent{} {\bf Notations:} Let $G(V,E)$ be a simple, finite, undirected graph on $n$ vertices. The vertex set of $G$ is denoted as $V(G)$ and the edge set of $G$ is denoted as $E(G)$. For any vertex $v \in V(G)$ let $N_G(v)=\{w \in V(G) \ | \ (v,w) \in E(G)\}$ be the set of neighbors of $v$. For each $S \subseteq V(G)$ let $G[S]$ denote the subgraph of $G$ induced by the vertices in $S$. In this paper we shall use the notation $G\setminus S$ to denote $G[V \setminus S]$. The \textit{girth} of a graph is the length of a shortest cycle in the graph.

Let $G'$ be a graph such that $V(G') = V(G)$. Then $G'$ is a {\it super graph} of $G$ if $E(G) \subseteq E(G')$. We define the {\it intersection} of two graphs as follows: if $G_1$ and $G_2$ are two graphs such that $V(G_1) = V(G_2)$, then the intersection of $G_1$ and $G_2$ denoted as $G = G_1 \cap G_2$ is a graph with $V(G) = V(G_1) = V(G_2)$ and $E(G) = E(G_1) \cap E(G_2)$.

\subsection{Boxicity and Cubicity}
A $d$-dimensional box is a Cartesian product $R_1\times R_2\times\cdots\times R_d$ where each $R_i$ is a closed interval of the form $[a_i,b_i]$ on the real line. A $k$-box representation of a graph $G$ is a mapping of the vertices of $G$ to $k$-boxes such that two vertices in $G$ are adjacent if and only if their corresponding $k$-boxes have a non-empty intersection. The \emph{boxicity} of a graph $G$, denoted as $\boxi(G)$, is the minimum integer $k$ such that $G$ can be represented as the intersection graph of $k$-dimensional boxes. Clearly, graphs with boxicity 1 are precisely the \emph{interval graphs}.


A $d$-dimensional cube is a Cartesian product $R_1\times R_2\times\cdots\times R_d$ where each $R_i$ is a closed interval of the form $[a_i,a_i+1]$ on the real line. A $k$-cube representation of a graph $G$ is a mapping of the vertices of $G$ to $k$-cubes such that two vertices in $G$ are adjacent if and only if their corresponding $k$-cubes have a non-empty intersection. The \emph{cubicity} of $G$ is the minimum integer $k$ such that $G$ has a $k$-cube representation. Clearly, graphs with cubicity 1 are precisely the \emph{unit interval graphs}.

Let $G$ be a graph. Let $I_1, I_2, \ldots, I_k$ be $k$ interval (unit interval) graphs such that $G = I_1 \cap I_2 \cap \cdots \cap I_k$. Then $I_1,
I_2, \ldots, I_k$ is called an {\it interval graph representation} ({\it unit interval graph representation}) of $G$. The following equivalence is well known.

\begin{lemma}\textbf{\textnormal{(Roberts\cite{recentProgressesInCombRoberts})}}
The minimum $k$ such that there exists an interval graph representation $($unit interval graph representation$)$ of $G$ using $k$ interval graphs $($unit interval graphs$)$ $I_1, I_2, \ldots, I_k$ is the same as $\boxi(G)$ $(\cub(G))$.
\end{lemma}

\begin{fact} \textnormal{(Roberts \cite{recentProgressesInCombRoberts})}\label{fact:cub_sum}
 If $G=G_1\cap G_2\cap\cdots\cap G_r$ then $\cub(G) \le \sum_{i=1}^r \cub(G_i)$.
\end{fact}

The concept of boxicity and cubicity was introduced by F. S. Roberts \cite{recentProgressesInCombRoberts} in 1969. Boxicity finds applications in fields such as ecology and operations research: It is used as a measure of the complexity of ecological \cite{mathematicalModelRoberts} and social \cite{spheresCubesFreeman} networks and has applications in fleet maintenance \cite{fleetOpsutRoberts}. Boxicity and cubicity has been investigated for various classes of graphs \cite{computingBoxicityCozzensRoberts,phdThesisScheinerman,intervalRepPlanarGraphsThomassen,cubicityMannino} and has been related with other parameters such as treewidth \cite{boxicityTreewidthSunilNaveen} and vertex cover \cite{boxicityVertexcoverSunil}. Computing the boxicity of a graph was shown to be NP-hard by Cozzens \cite{computingBoxicityCozzensRoberts}. This was later strengthened by Yannakakis \cite{complexityPartialOrderDimnYannakakis}, and finally by Kratochv\'{\i}l \cite{specialPlanarSatisfiabilityProbNPKratochvil} who showed that deciding whether boxicity of a graph is at most $2$ itself is NP-complete. Boxicity has been generalized in several ways like rectangle number \cite{rectNumHypercubesChangeWest}, poset boxicity \cite{posetBoxicityTrotterWest}, grid dimension \cite{gridIntBoxicityBellantoniEtal}, circular dimension \cite{circDimensionFeinberg}, boxicity of digraphs \cite{intNumBoxDigraphs} etc. Recently Chandran et al. \cite{boxicityMaxDegreeSunilNaveenFrancis} showed that for any graph $G$, $\boxi(G) \le 2\chi(G^2)$ where $G^2$ is the square of graph $G$ and $\chi(G)$ is the chromatic number of the graph. From this they inferred that $\boxi(G) \le 2\Delta^2$, where $\Delta$ is the maximum degree of $G$. Very recently this result was improved by Esperet \cite{boxicityEsperet}, who showed that $\boxi(G) \le \Delta^2+2$. Let $n$ be the number of vertices in $G$. In \cite{geometricSunil} Chandran et al. have shown that for any graph $G$, $\boxi(G) \le \left \lceil(\Delta+2)\log_2 n \right \rceil$. In \cite{bandwidth} they have shown that for any graph $G$, $\cub(G) \le \left \lceil 4(\Delta+1)\log_2 n \right \rceil$.

\subsection{Chordal Graph and Chordal Dimension}
An undirected graph is said to be \textit{chordal} if every cycle of length four or more contains a chord i.e. an edge joining two nonconsecutive vertices in the cycle. The \textit{chordal dimension} of a graph $G$ denoted as $\cd(G)$, is the minimum integer $k$ such that $G$ can be represented as the intersection graph of $k$ chordal graphs. Scheinerman and Mckee \cite{chordalityMcKeeScheinerman} have shown that for any graph $G$, $\cd(G) \le \chi(G)$ and also $\cd(G) \le treewidth(G)$ where $\chi(G)$ is the chromatic number of $G$. Since any interval graph is a chordal graph we have the following observation:

\begin{observation}
 For any graph $G$, $\cd(G) \le \boxi(G) \le \cub(G)$.
\end{observation}

\subsection{Claw Number}
Let $S(k)$ denote a star graph on $k+1$ vertices. (Note that $S(k)$ is the complete bipartite graph $K_{1,k}$). The \textit{center} of a star is that vertex which is connected to all other vertices in the star. An induced $S(3)$ in a graph is usually known as a \textit{claw}.
\begin{definition}
 The \emph{claw number} of a graph $G$ is the largest positive integer $k$ such that $S(k)$ is an induced subgraph of $G$ and is denoted as $\psi(G)$.
\end{definition}
Recently Adiga et al. \cite{cubInterval} have given an almost tight bound for the cubicity of interval graphs in terms of its claw number.
\begin{theorem}\textnormal{(Adiga et al. \cite{cubInterval})}\label{theo:claw}
 If $G$ is an interval graph with claw number $\psi(G)$ then $\ceil{\log_2\claw} \le \cub(G) \le \ceil{\log_2\claw}+2$.
\end{theorem}
\subsection{AT-free graphs}
An independent set of three vertices is called an \textit{asteroidal triple} if between every pair of vertices there is a path which avoids the neighbourhood of the third. A graph is called \textit{asteroidal triple free} (AT-free for short) if it does not contain an asteroidal triple. They form a large class of graphs since they contain interval, permutation, trapezoid, cocomparability and many other graph classes. Corneil, Olariu and Stewart have studied many structural and algorithmic properties of AT-free graphs in \cite{atfreeCorneil,dominatingCorneil}.

A graph is called \textit{claw-free AT-free graph} if it is AT-free and does not contain $K_{1,3}$ (i.e. $S(3)$, the claw) as an induced subgraph. Kloks et al. \cite{atfreeKloks} have given a characterization of claw-free AT-free graphs.

\subsection{Our Results}
In this paper we will show that
\begin{enumerate}
 \item If $G$ is an AT-free graph with chromatic number $\chi(G)$ then $\boxi(G) \le \chi(G)$ and this bound is tight.
 \item If $G$ is a claw-free AT-free graph with chromatic number $\chi(G)$ then $\boxi(G)=\cub(G) \le \chi(G)$ and this bound is tight.
 \item If $G$ is an AT-free graph having girth at least $5$ then $\boxi(G) \le 2$ and this bound is tight. We also show that $\cub(G) \leq 2\ceil{\log_2 \claw} +4$.
 \item If $G$ is an AT-free graph with chromatic number $\chi(G)$ and claw number $\claw$ then $\cub(G) \le \boxi(G)(\ceil{\log_2 \claw} +2)$ $\le \chi(G)(\ceil{\log_2 \claw} +2)$.
\end{enumerate}

\noindent{}\textbf{Remark on previous approach to boxicity and cubicity of AT-free graphs:}\\
In \cite{boxicityTreewidthSunilNaveen} it has been shown that for any graph $G$, $\boxi(G) \leq treewidth(G)+2$. It has also been shown that if $G$ is an AT-free graph then $treewidth(G) \leq 3\Delta-2$, hence $\boxi(G) \leq 3\Delta$ where $\Delta$ is the maximum degree of $G$. But the result shown in this paper is much stronger. (Recall that $\chi(G) \leq \Delta +1$ for any graph, but in general $\chi(G)$ can be much smaller.)

In \cite{bandwidth} Chandran et al. have studied the relationship between cubicity and bandwidth of a graph. As a corollary they have also shown that if $G$ is an AT-free graph then $\cub(G) \le 3\Delta-1$ since for AT-free graphs bandwidth is at most $3\Delta-2$. Using the technique of \cite{bandwidth}, this upper bound cannot be improved much since $\ceil{\frac{\Delta}{2}}$ is a lower bound for bandwidth of any graph. In this paper we show that for any AT-free graph $G$, $\cub(G) \le \boxi(G)(\ceil{\log_2 \claw} +2)$ $\le \chi(G)(\ceil{\log_2 \claw} +2)$. Clearly this result can be much stronger than that of \cite{bandwidth} in some cases.

\section{Upper bound on boxicity of AT-free graphs and cubicity of claw-free AT-free graphs}
In this section we will show an upper bound on boxicity of AT-free graphs and cubicity of claw-free AT-free graphs. A \textit{triangulation} of a graph $G$ is a chordal graph $H$ on the same vertex set that contains $G$ as a subgraph i.e. $V(G)=V(H)$ and $E(G) \subseteq E(H)$. $H$ is said to be a \textit{minimal triangulation} of $G$ if there exists no other chordal graph $H'$ on the same vertex set as $G$ and $H$ such that $E(G) \subseteq E(H') \subset E(H)$. M\"{o}hring studied minimal triangulation of AT-free graphs in \cite{triangulationMohring}. Parra and Schefller have shown relations between minimal separators of a graph and its minimal triangulations in \cite{characterizationChordalParra}.

From the definition of chordal dimension and boxicity we know that for any graph $G$, $\cd(G) \le \boxi(G)$. Now we will show that when $G$ is an AT-free graph, $\boxi(G) \le \cd(G)$. For this we need the following theorem:

\begin{theorem}\textnormal{(M\"{o}hring \cite{triangulationMohring})}\label{mohring}
 If $G$ is an AT-free graph then every minimal triangulation of $G$ is an interval graph.
\end{theorem}
Let $\cd(G)=k$ and $G=\bigcap_{i=1}^k G_i$ where $G_i$ is a chordal graph for $1 \le i \le k$. It is easy to see that if we replace each $G_i$ by another chordal graph $G_i'$ such that $V(G_i)=V(G_i')$ and $E(G) \subseteq E(G_i') \subseteq E(G_i)$, we still will have $G=\bigcap_{i=1}^k G_i'$. It follows that there exists $G_1',G_2',\cdots,G_k'$ such that $G=\bigcap_{i=1}^k G_i'$ where each $G_i'$ is a minimal triangulation of $G$. By Theorem \ref{mohring}, $G_i'$ for $1 \le i \le k$ is an interval graph. It follows that $\boxi(G) \le k= \cd(G)$. Thus we have the following Observation:

\begin{observation}\label{obs:box}
 If $G$ is an AT-free graph then $\cd(G)=\boxi(G)$.
\end{observation}

Scheinerman and Mckee have shown the following upper bound on chordal dimension of a graph $G$ in terms of its chromatic number $\chi(G)$.

\begin{theorem}\textnormal{(Scheinerman and Mckee \cite{chordalityMcKeeScheinerman})}\label{chordality}
 For any graph $G$ with chromatic number $\chi(G)$, $\cd(G) \le \chi(G)$.
\end{theorem}

Combining Observation \ref{obs:box} and Theorem \ref{chordality} we get the following upper bound on boxicity of AT-free graphs:

\begin{theorem}\label{box_general}
 If $G$ is an AT-free graph with chromatic number $\chi(G)$ then $\boxi(G) \le \chi(G)$.
\end{theorem}
 In general $\chi(G) \leq d+1$, where $d$ is the degeneracy of the graph. It follows that $box(G) \leq d+1$. Though it is known \cite{boxicityMaxDegreeSunilNaveenFrancis} that $\boxi(G) \leq 2\chi(G^2)$ for any graph $G$, $\boxi(G)$ need not always be less than equal to $\chi(G)$: For example it is shown in \cite{boxicityVertexcoverSunil} that there exists bipartite graphs with boxicity $\frac{n}{4}$. It is also shown in \cite{lowerBoundBoxicity} that almost all balanced bipartite graphs (with respect to a suitable probability distribution) have boxicity $\Omega(\frac{n}{\log{n}})$.

\begin{theorem}\textnormal{(Parra and Scheffler \cite{characterizationChordalParra})}\label{parra}
 A graph $G$ is claw-free AT-free if and only if every minimal triangulation of $G$ is a unit interval graph.
\end{theorem}

\noindent{}By a similar argument given for Observation \ref{obs:box}, we get the following:
\begin{observation}\label{obs:claw_box}
 If $G$ is a claw-free AT-free graph then $\cd(G)=\cub(G)$.
\end{observation}

\noindent{}Thus if $G$ is a claw-free AT-free graph we have $\cd(G)=\boxi(G)=\cub(G)$. Combining Theorem \ref{chordality} and Observation \ref{obs:claw_box} we get the following upper bound on cubicity of claw-free AT-free graphs:

\begin{theorem}\label{cubicity_general}
 If $G$ is a claw-free AT-free graph with chromatic number $\chi(G)$ then $\cub(G) \le \chi(G)$.
\end{theorem}

\subsection{Tightness of Theorem \ref{box_general} and Theorem \ref{cubicity_general}}
Let $G$ be a complete $k$-partite graph on $n$ vertices (We will assume that $n$ is multiple of $k$ and $n >k$). It is easy to see that this is an AT-free graph. Since the chromatic number of this graph is $k$, we have $\boxi(G) \leq k$ by Theorem \ref{box_general}. But it was shown by Roberts \cite{recentProgressesInCombRoberts} that $\boxi(G) = k$. So the upper bound for boxicity given in Theorem \ref{box_general} is tight for complete $k$-partite graphs.

Let $G=\overline{(\frac{n}{2})K_2}$, the complement of the perfect matching on $n$ vertices (We will assume that $n$ is even and $n >3$). It is easy to see that this is a claw-free AT-free graph. Since the chromatic number of this graph is $\frac{n}{2}$, we have $\cub(G) \leq \frac{n}{2}$ by Theorem \ref{cubicity_general}. But it was shown by Roberts \cite{recentProgressesInCombRoberts} that $\cub(G) = \frac{n}{2}$. So the upper bound for
cubicity given in Theorem \ref{cubicity_general} is tight for $\overline{(\frac{n}{2})K_2}$.

\section{Upper bound on boxicity of AT-free graphs having girth at least 5}
In this section we will show an upper bound on boxicity of AT-free graphs having girth at least $5$. Let $G$ be an AT-free graph having girth at least $5$. Since an induced cycle of length at least $6$ contains an AT, $G$ is either acyclic or all induced cycles of $G$ are of length exactly $5$. Recall that diameter of a graph is the maximum of distance(u,v) over all pairs of vertices $u,v \in V(G)$. A set of vertices $S$ of a graph $G$ is said to be \textit{dominating} if every vertex in $V(G) \setminus S$ is adjacent to some vertex in $S$. A path joining vertices $x$ and $y$ is said to be a $x$-$y$ path. A pair of vertices $x,y$ is said to be a \textit{dominating pair} if all $x$-$y$ paths in $G$ are dominating sets. Corneil, Olariu and Stewart have shown the following fundamental property of AT-free graphs:

\begin{theorem}\textbf{\textnormal{(Corneil et al.\cite{atfreeCorneil})}}
 Every connected AT-free graph contains a dominating pair.
\end{theorem}

\noindent{} They have also proved the following theorem which we shall use to show the upper bound on boxicity:
\begin{theorem}\textbf{\textnormal{(Corneil et al.\cite{atfreeCorneil})}}\label{theo:dompair}
 In every connected AT-free graph there exists dominating pair $x,y$ such that $distance(x,y)=diameter(G)$.
\end{theorem}

\noindent{} Let $x,y$ be a dominating pair in $G$ and let $P$ be a shortest $x$-$y$ path of length equal to the diameter of $G$. Let $d$ be the diameter of $G$ and $V(P)=\{u_1,u_2,\cdots,u_d\}$ where $x=u_1$ and $y=u_d$. Let $V(\comp{P})=V(G) \setminus V(P)$.

\begin{lemma}\label{lem:single}
 For each vertex  $v \in V(\comp{P})$, $|N_G(v) \cap V(P)| =1$.
\end{lemma}
\begin{proof}
Since $x,y$ is a dominating pair and $P$ is a $x$-$y$ path, $V(P)$ is a dominating set. Hence for every vertex $v \in V(\comp{P})$ we have $|N_G(v) \cap V(P)| \ge 1$. We will show that for each vertex  $v \in V(\comp{P})$, $|N_G(v) \cap V(P)| \le 1$. If possible let $w \in V(\comp{P})$ such that $|N_G(w) \cap V(P)| \ge 2$. Let $u_i,u_j \in N_G(w) \cap V(P)$ be such that $1 \leq i<j \leq d$ and for all $k$, $i < k<j$, $u_k \notin N_G(w)$. We consider the following cases \\

\noindent{}\textbf{Case 1:}{\it When $j \le i+2$.} If $j=i+1$ then $u_i$-$w$-$u_j$-$u_i$ forms an induced cycle of length $3$ in $G$, a contradiction. Similarly if $j=i+2$ then $u_i$-$w$-$u_j$-$u_{j-1}$-$u_i$ forms an induced cycle of length $4$ in $G$, a contradiction.

\noindent{}\textbf{Case 2:}{\it When $j \ge i+3$.} Let $P_1$ denote the path $u_1$-$u_2$-$\cdots$-$u_i$ and $P_2$ denote the path $u_j$-$u_{j+1}$-$\cdots$-$u_d$. Clearly $P_1wP_2$ forms a $x$-$y$ path say $P'$ in $G$. Now $|V(P')| = i+1+(d-j+1)$. If $j \ge i+3$ then $|V(P')| \le d-1$. Recall that $P$ is a shortest $x$-$y$ path. But $P'$ is a shorter $x$-$y$ path than $P$, a contradiction.\\

\noindent{} Therefore for each vertex  $v \in V(\comp{P})$, $|N_G(v) \cap V(P)| = 1$.
\qed
\end{proof}

Let $S_i=\{v \ | \ v \in V(\comp{P}) \ and \ u_i \in N_G(v)\}$ for $1
\le i \le d$. From Lemma \ref{lem:single}, it follows that $S_1,S_2,\cdots,S_d$ is a partition of the vertex set $V(\comp{P})$. In other words,

\begin{observation}\label{obs:nonadj}
 Let $u \in V(P)$ and $v \in V(\comp{P})$. Suppose $u=u_i$ and $v \in S_k$ where $1 \leq i,k \leq d$. Then $(u,v) \notin E(G)$ if and only if $i \neq k$.
\end{observation}

\begin{lemma}\label{Lem:adj}
Let $v \in S_i$.
\begin{enumerate}
 \item $|N_G(v) \cap S_i|=0$ where $1 \le i \le d$.
 \item $|N_G(v) \cap S_{i+1}|=0$ where $1 \leq i \leq d-1$.
 \item $|N_G(v) \cap S_{i+2}|\le 1$ where $1 \leq i \leq d-2$.
 \item $|N_G(v) \cap S_{j}|=0$ where $i+3 \leq j \leq d$ and $i \geq 1$.
\end{enumerate}
\end{lemma}

 \noindent{}\textbf{Proof(1):} If possible let $w \in S_i$ such that $(v,w) \in E(G)$. Now $v$-$u_i$-$w$-$v$ forms an induced cycle of length $3$ in $G$, a contradiction.\\

\noindent{}\textbf{Proof(2):} If possible let $w \in S_{i+1}$ such that $(v,w) \in E(G)$. Then $u_i$-$v$-$w$-$u_{i+1}$-$u_i$ forms a cycle of length $4$ in $G$, a contradiction.\\

\noindent{}\textbf{Proof(3):} If possible let $u,w \in S_{i+2}$ such that $(v,u) \in E(G)$ and $(v,w) \in E(G)$. Then $v$-$w$-$u_{i+2}$-$u$-$v$ forms a cycle of length $4$ in $G$, a contradiction.\\

\noindent{}\textbf{Proof(4):} If possible let $w \in S_j$ such that $(v,w) \in E(G)$. Since $(v,u_i) \in E(G)$ according to Lemma \ref{lem:single} we have $(v,u_k) \notin E(G)$ for all $k \ne i$. Similarly since $(w,u_j) \in E(G)$ we have $(w,u_k) \notin E(G)$ for all $k \ne j$. Since $j \ge i+3$, $u_{i}$-$v$-$w$-$u_{j}$-$u_{j-1}$-$u_{j-2}$-$\cdots$-$u_i$ forms an induced cycle of length at least $6$ in $G$. But $G$ is an AT-free graph, a contradiction.


\noindent{} From Lemma \ref{Lem:adj} we have the following observation:

\begin{observation}\label{obs:adj}
 If $u,v \in V(\comp{P})$, $u \in S_i$, $v \in S_j$ and $(u,v) \in E(G)$ then $|j-i|=2$.
\end{observation}

\begin{lemma}\label{Lem:special}
 Let $u \in S_i$ and $v \in S_{i+2}$ where $1 \leq i \leq d-2$. If $(u,v) \in E(G)$ then for any $p \in S_i \setminus \{u\}$, $q \in S_{i+2} \setminus \{v\}$ we have $(p,q) \notin E(G)$.
\end{lemma}

\begin{proof}
Suppose not. Let $p \in S_i \setminus \{u\}$ and $q \in S_{i+2} \setminus \{v\}$ such that $(p,q) \in E(G)$. Since $u,p \in S_i$, according to Lemma \ref{Lem:adj} part (1), $(u,p) \notin E(G)$. Similarly $(v,q) \notin E(G)$. According to Lemma \ref{lem:single}, $(u,u_{i+2}) \notin E(G)$, $(p,u_{i+2}) \notin E(G)$, $(q,u_i) \notin E(G)$ and $(v,u_i) \notin E(G)$. Also we have $(u,q) \notin E(G)$ and $(v,p) \notin E(G)$ by Lemma \ref{Lem:adj} part (3). Moreover $(u_i,u_{i+2}) \notin E(G)$ since $P$ is a shortest $x$-$y$ path. Therefore $u$-$u_i$-$p$-$q$-$u_{i+2}$-$v$-$u$ forms an induced cycle of length $6$ in $G$ and hence $\{u,p,u_{i+2}\}$ forms an AT in $G$, a contradiction.
 \qed
\end{proof}


A vertex $v \in V(\comp{P})$ is said to be \textit{non-pendant} if $N_G(v) \cap V(\comp{P}) \ne \emptyset$. Note that if $N_G(v) \cap V(\comp{P})=\emptyset$ then $v$ has to be a pendant vertex by Lemma \ref{lem:single}. 

\begin{lemma}
 $S_i$ can contain at most $2$ non-pendant vertices for $1 \leq i \leq d$.
\end{lemma}

\begin{proof}
If $v \in S_i$ is non-pendant then according to Observation \ref{obs:adj}, either $N_G(v) \cap S_{i-2} \neq \emptyset$ or $N_G(v) \cap S_{i+2} \neq \emptyset$. By Lemma \ref{Lem:adj} part (3) and Lemma \ref{Lem:special}, at most one vertex in $S_i$ can be connected to some vertex in $S_{i+2}$. Similarly at most one vertex in $S_i$ can be connected to some vertex in $S_{i-2}$. Therefore $S_i$ can contain at most $2$ non-pendant vertices.
\qed
\end{proof}

\begin{observation}\label{obs:special2}
 If $S_i$ contains two non-pendant vertices say $u,v$ then one of the following statements is true $($by Lemma \ref{Lem:adj} part $(3)$ and Lemma \ref{Lem:special}$)$
\begin{enumerate}
 \item $N_G(u) \cap S_{i-2} = \emptyset$ and $N_G(v) \cap S_{i+2} = \emptyset$.
 \item $N_G(u) \cap S_{i+2} = \emptyset$ and $N_G(v) \cap S_{i-2} = \emptyset$.
\end{enumerate}
\end{observation}

\subsection{Interval Graph Construction}
We shall construct two interval graphs $I_1$ and $I_2$ such that $G=I_1 \cap I_2$. In the interval graph $I_j$ where $j=1,2$, let $l_j(u)$ and $r_j(u)$ denote the left and right endpoint of the interval corresponding to vertex $u \in V(G)$ respectively. Let $S$ be the set of non-pendant vertices in $V(\comp{P})$. To construct $I_1$ we map each vertex $v \in V(G)$ to an interval on the real line by the mapping:\\

\vspace{-0.6cm}
\begin{eqnarray*}
    g_1(v) &=& [i,i+1] \ \ \ \ \ \ \ \ \ \ \ \ \ \ \ \ \ \ if \ v \in V(P) \ and \ v=u_i \ for \ 1 \le i \le d\\
    &=& [i+\frac{2j-1}{2n},i+\frac{2j}{2n}] \ \ \ \ \ if \ v \in S_i \setminus S, \ 1 \le i \le d \ and \ 1 \le j \le |S_i|\\
    &=& [i-\frac{1}{2},i+\frac{3}{2}] \ \ \ \ \ \ \ \ \ \ \ \ if \ v \in S_i \cap S, \ N_G(v) \cap S_{i-2} \neq \emptyset \ and  \ N_G(v) \cap S_{i+2} \neq \emptyset\\
    &=& [i+1,i+\frac{3}{2}] \ \ \ \ \ \ \ \ \ \ \ \ \ if \ v \in S_i \cap S, \ N_G(v) \cap S_{i-2} = \emptyset \ and  \ N_G(v) \cap S_{i+2} \neq \emptyset\\
    &=& [i-\frac{1}{2},i] \ \ \ \ \ \ \ \ \ \ \ \ \ \ \ \ \ \ if \ v \in S_i \cap S, \ N_G(v) \cap S_{i-2} \neq \emptyset \ and  \ N_G(v) \cap S_{i+2} = \emptyset
\end{eqnarray*}

\begin{lemma}\label{Lem:sup1}
 $I_1$ is a supergraph of $G$.
\end{lemma}
\begin{proof}
Let $(u,v) \in E(G)$. We shall show that $g_1(u) \cap g_1(v) \neq \emptyset$. We consider the following cases:

\noindent{}\textbf{Case 1:} {\it When either $u \in V(P)$ or $v \in V(P)$.} Without loss of generality we can assume that $u \in V(P)$. Let $u=u_i$ where $1 \leq i \leq d$. If $v \in V(P)$ then either $v=u_{i-1}$ or $v=u_{i+1}$. Now if $v=u_{i-1}$ then $i \in g_1(u) \cap g_1(v)$. On the other hand if $v=u_{i+1}$ then $i+1 \in g_1(u) \cap g_1(v)$. If $v \in \comp{P}$ then according to Observation \ref{obs:nonadj}, $v \in S_i$. Now if $v \in S_i \setminus S$ then $i=l_1(u) < l_1(v) < r_1(u)=i+1$ and hence $g_1(u) \cap g_1(v) \neq \emptyset$. If $v \in S_i \cap S$ then we consider the following cases: If $N_G(v) \cap S_{i-2} \neq \emptyset$ then $i \in g_1(u) \cap g_1(v)$. Otherwise if $N_G(v) \cap S_{i+2} \neq \emptyset$ then $i+1 \in g_1(u) \cap g_1(v)$.\\

\noindent{}\textbf{Case 2:} {\it When $u,v \in V(\comp{P})$.} By definition of non-pendant vertices $u,v \in S$. Let $u \in S_i$. According to Observation \ref{obs:adj}, either $v \in S_{i-2}$ or $v \in S_{i+2}$. If $v \in S_{i-2}$ then $l_1(u)=r_1(v)=i-\frac{1}{2}$. Otherwise if $v \in S_{i+2}$ then $r_1(u)=l_1(v)=i+\frac{3}{2}$. In both cases we have $g_1(u) \cap g_1(v) \neq \emptyset$.
 \qed
\end{proof}

To construct $I_2$ we map each vertex $v \in V(G)$ to an interval on the real line by the mapping:\\

\vspace{-.6cm}
\begin{eqnarray*}
    g_2(v) &=& [1,2] \ \ \ \ \ \ \ \ \ \ if \ v \in V(P) \ , \ v=u_i \ and \ i \ mod \ 2 =1\\
    &=& [2,3] \ \ \ \ \ \ \ \ \ \ if \ v \in V(P) \ , \ v=u_i \ and \ i \ mod \ 2 =0\\
    &=& [\frac{5}{4},\frac{7}{4}] \ \ \ \ \ \ \ \ \ if \ v \in S_i \setminus S \ and \ i \ mod \ 2 =1\\
    &=& [\frac{9}{4},\frac{11}{4}] \ \ \ \ \ \ \ if \ v \in S_i \setminus S \ and \ i \ mod \ 2 =0\\
    &=& [0,1] \ \ \ \ \ \ \ \ \ \ if \ v \in S_i \cap S \ and \ i \ mod \ 2 =1\\
    &=& [3,4] \ \ \ \ \ \ \ \ \ \ if \ v \in S_i \cap S \ and \ i \ mod \ 2 =0
\end{eqnarray*}

\begin{lemma}\label{Lem:sup2}
 $I_2$ is a supergraph of $G$.
\end{lemma}

\begin{proof}
Let $(u,v) \in E(G)$. We shall show that $g_2(u) \cap g_2(v) \neq \emptyset$. We consider the following cases:

\noindent{}\textbf{Case 1:} {\it When either $u \in V(P)$ or $v \in V(P)$.} Without loss of generality we can assume that $u \in V(P)$. Let $u=u_i$ where $1 \leq i \leq d$. If $v \in V(P)$ then $2 \in g_2(u) \cap g_2(v)$. If $v \in V(\comp{P})$ then according to Observation \ref{obs:nonadj}, $v \in S_i$. Now if $v \in S_i \setminus S$ then $l_2(u) < l_2(v) <r_2(v) < r_2(u)$ and hence $g_2(u) \cap g_2(v) \neq \emptyset$. If $v \in S_i \cap S$ we consider the following cases: If $i \ mod \ 2=1$ then $l_2(u)=r_2(v)=1$. On the other hand if $i \ mod \ 2=0$ then $r_2(u)=l_2(v)=3$. In both cases we have $g_2(u) \cap g_2(v) \neq \emptyset$.\\

\noindent{}\textbf{Case 2:} {\it When $u,v \in V(\comp{P})$.} By definition of non-pendant vertices $u,v \in S$. Let $u \in S_i$ and $v \in S_j$ where $1 \leq i,j \leq d$. According to Observation \ref{obs:adj}, $|i-j|=2$ which implies that $i = j \ mod \ 2$. Hence $g_2(u)=g_2(v)$ and thus $g_2(u) \cap g_2(v) \neq \emptyset$.
 \qed
\end{proof}

\begin{lemma}\label{Lem:misssing_edge}
 For any $(u,v) \notin E(G)$ either $(u,v) \notin E(I_1)$ or $(u,v) \notin E(I_2)$.
\end{lemma}

\begin{proof}
Let $(u,v) \notin E(G)$. We consider the following cases:

\noindent{}\textbf{Case 1:} {\it When $u,v \in V(P)$.} Let $u=u_i$ and $v =u_j$ where $1 \leq i,j \leq d$. Since $(u,v) \notin E(G)$ we have $|j-i|\geq 2$. Therefore $|l_1(u)-l_1(v)| \geq 2$. Since in $I_1$, the intervals corresponding to vertices in $V(P)$ are of length $1$ we have $g_1(u) \cap g_1(v) = \emptyset$ and hence $(u,v) \notin E(I_1)$.\\

 \noindent{}\textbf{Case 2:} {\it When $u \in V(P)$ and $v \in V(\comp{P})$.} Let $u=u_i$ where $1 \leq i \leq d$. According to Observation \ref{obs:nonadj}, $v \in S_k$ where $k \neq i$. Now if $v \in S_k \setminus S$ then $k<l_1(v) < r_1(v) <k+1$. Since $g_1(u)=[i,i+1]$ and $i \neq k$ we have $g_1(u) \cap g_1(v) = \emptyset$. Hence $(u,v) \notin E(I_1)$.\\ 
When $v \in S_k \cap S$ we consider the following cases:

\noindent{}\textbf{Subcase 2.1:} {\it When $|k-i| \geq 2$.} Now $g_1(u)=[i,i+1]$ and $k-\frac{1}{2} \leq l_1(v) < r_1(v) \leq k+\frac{3}{2}$. If $i \leq k-2$ then $r_1(u) \leq k-1 <k-\frac{1}{2} \leq l_1(v)$ and hence $g_1(u) \cap g_1(v) = \emptyset$. If $i \geq k+2$ then $l_1(u) \geq k+2 >k+\frac{3}{2} \geq r_1(v)$ and hence $g_1(u) \cap g_1(v) = \emptyset$. Therefore $(u,v) \notin E(I_1)$.

\noindent{}\textbf{Subcase 2.2:} {\it When $|k-i| \leq 1$.} Since $k \neq i$ we have $k \ mod \ 2 \neq i \ mod \ 2$. If $i \ mod \ 2=0$ then $g_2(u)=[2,3]$ and $g_2(v)=[0,1]$. Hence $g_2(u) \cap g_2(v) = \emptyset$. If $i \ mod \ 2=1$ then $g_2(u)=[1,2]$ and $g_2(v)=[3,4]$. Hence $g_2(u) \cap g_2(v) = \emptyset$. In both cases we have $(u,v) \notin E(I_2)$.\\

\noindent{}\textbf{Case 3:} {\it When $u,v \in V(\comp{P})$.} We consider the following cases:\\
\noindent{}\textbf{Subcase 3.1:} {\it When $u,v \in S$.} Let $u \in S_i$ and $v \in S_j$. If $i=j$ then according to Observation \ref{obs:special2}, either $N_G(u) \cap S_{i-2} = \emptyset$ and $N_G(v) \cap S_{i+2} = \emptyset$ OR $N_G(v) \cap S_{i-2} = \emptyset$ and $N_G(u) \cap S_{i+2} = \emptyset$. If $N_G(u) \cap S_{i-2} = \emptyset$ and $N_G(v) \cap S_{i+2} = \emptyset$ then $r_1(v)=i<i+1=l_1(u)$. Hence $g_1(u) \cap g_1(v) = \emptyset$. If $N_G(v) \cap S_{i-2} = \emptyset$ and $N_G(u) \cap S_{i+2} = \emptyset$ then $r_1(u)=i<i+1=l_1(v)$. Hence $g_1(u) \cap g_1(v) = \emptyset$.

\noindent{}If $i\neq j$ then we consider the following cases. Without loss of generality we can assume that $j>i$.

\noindent{}\textbf{Subcase 3.1.1:}{\it When $(j-i) \ mod \ 2 \ne 0$.} It is easy to see that $g_2(u) \cap g_2(v) =\emptyset$.

\noindent{}\textbf{Subcase 3.1.2:}{\it When $(j-i) \ mod \ 2 = 0$.} We consider the following cases:

\noindent{}\textbf{Subcase 3.1.2.1:}{\it When $j=i+2$.} We will show that either $N_G(u)\cap S_{i+2}=\emptyset$ or $N_G(v) \cap S_{i}=\emptyset$. If possible let $N_G(u)\cap S_{i+2} \neq \emptyset$ and $N_G(v) \cap S_{i} \neq \emptyset$. Let $p \in S_i$ and $q \in S_{i+2}$ be such that $(u,q) \in E(G)$ and $(v,p) \in E(G)$. Since $(u,v) \notin E(G)$ we have $u \neq p$ and $v \neq q$. But then we get a contradiction to Lemma \ref{Lem:special}. Therefore either $N_G(u)\cap S_{i+2}=\emptyset$ or $N_G(v) \cap S_{i}=\emptyset$. If $N_G(u)\cap S_{i+2}=\emptyset$ then $r_1(u)=i < i+\frac{3}{2}=j-\frac{1}{2} \leq l_1(v)$. Therefore $g_1(u) \cap g_1(v) =\emptyset$. On the other hand if $N_G(v) \cap S_{i}=\emptyset$ then $r_1(u) \leq i+\frac{3}{2} < j+1 =l_1(v)$. Therefore $g_1(u) \cap g_1(v) =\emptyset$.

\noindent{}\textbf{Subcase 3.1.2.2:}{\it When $j \geq i+4$.} Then $r_1(u) \leq i+\frac{3}{2} <(i+4)-\frac{1}{2} \leq j-\frac{1}{2} \leq l_1(v)$. Therefore $g_1(u) \cap g_1(v) =\emptyset$. 

\noindent{}\textbf{Subcase 3.2:} {\it When $u \notin S$ and $v \notin S$.} According to the construction of $I_1$, it is easy to see that $\bigcup_{i=1}^d(S_i \setminus S)$ induces an independent set in $I_1$. Therefore $g_1(u) \cap g_1(v) =\emptyset$.

\noindent{}\textbf{Subcase 3.3:} {\it When $u \notin S$ and $v \in S$.} In $I_2$, $g_2(v)$ is either $[0,1]$ or $[3,4]$ and $g_2(u)$ is either $[\frac{5}{4},\frac{7}{4}]$ or $[\frac{9}{4},\frac{11}{4}]$. In all the four possible cases it is easy to see that $g_2(u) \cap g_2(v) =\emptyset$.

\qed
\end{proof}

\noindent{}Combining Lemma \ref{Lem:sup1},\ref{Lem:sup2} and \ref{Lem:misssing_edge} we have the following Theorem
\begin{theorem}\label{theo:box}
 If $G$ is an AT-free graph having girth at least $5$ then $\boxi(G) \le 2$.
\end{theorem}

\subsection{Tightness of Theorem \ref{theo:box}}
Let $G$ be a cycle of length $5$. It is easy to see that $G$ is an AT-free graph having girth at least $5$. According to Theorem \ref{theo:box}, $\boxi(G) \leq 2$. But clearly $\boxi(G)=2$, since $G$ is not an interval graph. Therefore the upper bound given by Theorem \ref{theo:box} is tight.

\section{Upper bound on cubicity of AT-free graphs}
In this section we will show an upper bound on cubicity of AT-free graphs in terms of its boxicity and claw number. This in turn will give an upper bound in terms of chromatic number and claw number. Let $G$ be an AT-free graph with chromatic number $\chi(G)$ and claw number $\claw$. We need some results shown by Parra and Scheffler \cite{characterizationChordalParra}.

For any graph $G(V,E)$ and for a given pair of nonadjacent vertices $a,b \in V$, a subset $S \subset V \setminus \{a,b\}$ is a $a$-$b$ \textit{vertex separator} ($a$-$b$ separator for short) if when $S$ is removed from $G$, $a$ and $b$ belong to different connected components of $G \setminus S$. $S$ is said to be a \textit{minimal} $a$-$b$ \textit{separator} if no proper subset of $S$ is an $a$-$b$ separator. A separator $S$ in $G$ is said to be a minimal separator of $G$ if there exists a pair of vertices $a,b \in V(G)$ such that $S$ is a minimal $a$-$b$ separator. It is well-known that a graph is chordal if and only if all its minimal separators induce cliques \cite{algGraphTheoryPerfectGraphsGolumbic}.

 Let $S$ and $T$ be two minimal separators of $G$. $S$ is said to \textit{cross} $T$ if there are two components $C,D$ of $G \setminus T$ such that $S$ intersects both $C$ and $D$. Parra and Schefller \cite{characterizationChordalParra} have shown that if $S$ crosses $T$ then $T$ crosses $S$ also. $S$ and $T$ are said to be \textit{parallel} if they do not cross each other. Let $\mathcal{S}_G$ denote the set of minimal separators in $G$. For $\mathcal{T}=\{S_1,S_2,\cdots,S_k\} \subseteq \mathcal{S}_G$, let $G_{\mathcal{T}}$ denote the graph obtained by making each separator $S_i$ for $1 \le i \le k$ a clique. The following Theorem is proved in \cite{characterizationChordalParra}.

\begin{theorem}\textnormal{(Parra and Scheffler \cite{characterizationChordalParra})}\label{separator}
\begin{enumerate}
 \item Let $\mathcal{T}=\{S_1,\cdots,S_k\}$ be a maximal set of pairwise parallel minimal separators in $G$. Then $H=G_{\mathcal{T}}$ is a minimal triangulation of $G$ and $\mathcal{S}_{H}=\mathcal{T}$.
 \item Let $H$ be a minimal triangulation of $G$. Then $\mathcal{S}_H$ is a maximal set of pairwise parallel minimal separators in $G$ and $H=G_{\mathcal{S}_{H}}$.
\end{enumerate}
\end{theorem}

%
%
 Let $T$ be a minimal separator of $G$. A component $C$ of $G \setminus T$ is called a \emph{full component} if every vertex in $T$ is adjacent to some vertex in $C$. The following property of minimal separator is shown in \cite{algGraphTheoryPerfectGraphsGolumbic}.
\begin{theorem}\textnormal{(Golumbic \cite{algGraphTheoryPerfectGraphsGolumbic})}\label{full_compo}
 A separator $T$ in graph $G$ is minimal if and only if there are at least two full components in $G \setminus T$.
\end{theorem}

\begin{lemma}\label{Lem:sep}
 Let $X$ be a minimal separator in a graph $G$ and $C$,$D$ be two full components in $G \setminus X$. Let $x \in X$, $c \in C$ and $d \in D$. Let $Y$ be another minimal separator of $G$ such that $c \in Y$ and $x,d \notin Y$. If $X$ is parallel to $Y$ then $x,d$ belongs to the same connected component in $G \setminus Y$.
\end{lemma}
\begin{proof}
 Suppose $x$ and $d$ lies in different connected components in $G \setminus Y$. Since $D$ is a full component in $G \setminus X$, there exists a $x$-$d$ path say $P$ in $G[D \cup \{x\}]$. Now according to assumption, $x$ and  $d$ lie in different components in $G \setminus Y$. Therefore $Y$ must contain at least one vertex from $P$. But since $x \notin Y$ and all the vertices in $P$ except $x$ lie in $D$ we have $Y \cap D \neq \emptyset$. Again $c \in Y \cap C$ and therefore $Y \cap C \neq \emptyset$. Hence $Y$ crosses $X$, a contradiction.
\qed
\end{proof}

\begin{lemma}\label{Lem:claw}
 If $G$ is an AT-free graph and $H$ is a minimal triangulation of $G$ with claw number $\psi(H)$ then $\psi(H) \le \claw$.
\end{lemma}
\begin{proof}
 Suppose $\psi(H)  > \claw$ and $\psi(H)=p$. An edge $(u,v) \in E(H)$ is said to be an old edge if $(u,v) \in E(G)$ and is said to be a new edge otherwise. Among all the claws of maximum size in $H$, let $U=\{s,x_1,x_2,\cdots,x_p\}$ induce the one with maximum number of old edges in it. Let $s$ be the center of the claw. Since $\psi(H) > \psi(G)$ at least one of the edges in $U$ has to be new. Without loss of generality let us assume that $(s,x_1)$ is a new edge. Let $\mathcal{T}=\{S_1,S_2,\cdots,S_k\}$ be the collection of minimal separators of $H$. From part (2) of Theorem \ref{separator}, $\mathcal{T}$ is a maximal set of pairwise parallel minimal separators of $G$ and $H=G_{\mathcal{T}}$. In other words if $(u,v) \in E(H) \setminus E(G)$ then there exists an $S_j \in \mathcal{T}$ such that both $u,v \in S_j$. Thus the vertices $s,x_1$ must belong to some minimal separator, say $X \in \mathcal{T}$ of $G$. Let $\mathcal{C}$ be the set of full components in $G \setminus X$. According to Theorem \ref{full_compo}, $|\mathcal{C}| \ge 2$. We consider the following two cases:\\

\noindent{}\textbf{Case 1:} {\it There exists a full component $C \in \mathcal{C}$ such that $C \cap \{x_2,x_3,\cdots,x_p\}=\emptyset$.} Since $C$ is a full component of $G \setminus X$ and $s \in X$ there is at least one vertex in $C$, say $a$ such that $(s,a) \in E(G)$. Since $E(G) \subseteq E(H)$ we have $(s,a) \in E(H)$. Note that $(a,x_i) \notin E(G)$ for $2 \le i \le p$ because $C \cap \{s,x_1,x_2,\cdots,x_p\}=\emptyset$ by assumption and $x_i \notin X$ for $2 \leq i \leq p$ since $x_1 \in X$ and $X$ induces a clique in $H$. Then it is easy to see that $\{s,a,x_2,\cdots,x_p\}$ forms a claw of size $p$ in $H$ having more old edges than in $U$ since $(s,x_1)$ is a new edge and $(s,a)$ is an old edge. But by assumption $U$ was a maximum sized claw having maximum number of old edges in it, a contradiction.\\

\noindent{}\textbf{Case 2:} {\it Every full component in $\mathcal{C}$ contains at least one $x_i$ where $ 2 \le i \le p$.} According to Theorem \ref{full_compo}, $|\mathcal{C}| \ge 2$ and hence there exists two full components $C,D \in \mathcal{C}$. Let $x_i \in C$ and $x_j \in D$ where $2 \le i < j \le p$. We will show that the triplet $\{x_1,x_i,x_j\}$ forms an AT in $G$, leading to a contradiction. Since $C$ is a full component of $G \setminus X$, $x_i \in C$ and $x_1 \in X$ there exists a $x_i$-$x_1$ path in $G[C\cup \{x_1\}]$ and this path does not intersect $N_G(x_j)$ since $x_j \in D$. Similarly since $D$ is a full component of $G \setminus X$, $x_j \in D$ and $x_1 \in X$ there exists a $x_j$-$x_1$ path in $G[D \cup \{x_1\}]$ which does not intersect $N_G(x_i)$. Now we want to show that there exists a $x_i$-$x_j$ path in $G$ which does not intersect $N_G(x_1)$. For that we need the following claim:\\

\noindent{}\textbf{Claim:}
\begin{enumerate}
 \item There exists a $x_i$-$s$ path in $G$ that does not intersect $N_G(x_1)$.
 \item There exists a $x_j$-$s$ path in $G$ that does not intersect $N_G(x_1)$.
\end{enumerate}
\begin{proof}
We prove only part (1) since the proof of part (2) is similar. Recall that $(s,x_1)$ is a new edge by assumption. Since $\{x_1,x_2,\cdots,x_p\}$ induce an independent set in $H$ and $E(G) \subseteq E(H)$ they induce an independent set in $G$ also. If $(s,x_i) \in E(G)$ we have a $x_i$-$s$ path in $G$ that does not intersect $N_G(x_1)$ since $(s,x_1) \notin E(G)$ and $(x_1,x_i) \notin E(G)$. Therefore we can assume that $(s,x_i) \notin E(G)$. Since $(s,x_i)$ is a new edge by theorem \ref{separator} there should be a minimal separator $Y \in \mathcal{T}$ such that $s,x_i \in Y$. Clearly $X \neq Y$ since $x_i \notin X$. According to Theorem \ref{separator}, $X$ and $Y$ are parallel and each separator in $\mathcal{T}$ induces a clique in $H$. Since $(x_i,x_1) \notin E(H)$, $(x_i,x_j) \notin E(H)$ and $x_i \in Y$ we have $x_1 \notin Y$ and $x_j \notin Y$. Therefore according to Lemma \ref{Lem:sep}, $x_1$ and $x_j$ must belong to the same connected component say $Q$ of $G \setminus Y$. Let $Q'$ be a full component of $G \setminus Y$ such that $Q' \neq Q$. Note that such a full component exists by Theorem \ref{full_compo}. Now $s$ and $x_i$ must be connected in $G$ to at least one vertex in $Q'$ and therefore there is a $x_i$-$s$ path in $G[Q' \cup \{x_i,s\}]$ which does not intersect $N_G(x_1)$. 
\qed
\end{proof}

Since $(s,x_1) \notin E(G)$ from the previous claim it is easy to see that there exists a $x_i$-$x_j$ path in $G$ which does not intersect $N_G(x_1)$. Therefore $\{x_1,x_i,x_j\}$ forms an asteroidal triple in $G$, a contradiction.
\qed
\end{proof}

\begin{theorem}\label{cub_general}
 If $G$ is an AT-free graph then $\cub(G) \leq \boxi(G)(\ceil{\log_2 \claw} +2) \leq \chi(G)(\ceil{\log_2 \claw} +2)$.
\end{theorem}
\begin{proof}
Let $\boxi(G)=k$ and $I_1,I_2,\cdots,I_k$ be interval graphs such that $G=\bigcap_{j=1}^{k}I_j$. It is easy to see that if we replace each $I_j$ by a chordal graph $I_j'$ such that $V(I_j)=V(I_j')$ and $E(G) \subseteq E(I_j') \subseteq E(I_j)$, we still have $G=\bigcap_{j=1}^{k}I_j'$. It follows that there exists chordal graphs $I_1',I_2',\cdots,I_k'$ such that $G=\bigcap_{j=1}^{k}I_j'$ where each $I_j'$ is a minimal triangulation of $G$. But by Theorem \ref{mohring} any minimal triangulation of an AT-free graph is an interval graph. It follows that $I_1',I_2',\cdots,I_k'$ are interval graphs. According to Lemma \ref{Lem:claw}, $\psi(I_j') \leq \claw$ for $1 \leq j \leq k$. Since $G=\bigcap_{j=1}^{k}I_j'$ we have $\cub(G) \leq \sum_{j=1}^k \cub(I_j')$ according to Fact \ref{fact:cub_sum}. By Theorem \ref{theo:claw},  $\cub(I_j') \leq \ceil{\log_2 \psi(I_j')} +2$  and by Lemma \ref{Lem:claw}, $\cub(I_j') \leq \ceil{\log_2 \psi(G)} +2$ for $1 \leq j \leq k$. It follows that $\cub(G) \leq k(\ceil{\log_2 \psi(G)} +2)=\boxi(G)(\ceil{\log_2 \psi(G)} +2)$. Therefore $\cub(G) \leq \boxi(G)(\ceil{\log_2 \claw} +2)$. By Theorem \ref{box_general}, we also have $\cub(G) \leq \chi(G)(\ceil{\log_2 \psi(G)} +2)$.
 \qed
\end{proof}

\noindent{}From Theorem \ref{theo:box} and \ref{cub_general} we get the following:
\begin{corollary}
 If $G$ is an AT-free graph having girth at least $5$ then $\cub(G) \leq 2\ceil{\log_2 \claw} +4$.
\end{corollary}

%
%
%
%
%
%
%

%
\end{document}